\documentclass[journal]{IEEEtran}

\usepackage[noadjust]{cite}

\usepackage{mathtools} % amsmath with fixes and additions
\usepackage{amsmath}
\usepackage{amssymb}
\usepackage{mathtools}
\usepackage{amsthm}
\usepackage{xcolor}
\usepackage{subcaption}
\usepackage{float}
\usepackage{bbm}
\usepackage[ruled,vlined]{algorithm2e}
\SetKwInput{KwInput}{Input}
\SetKwInput{KwOutput}{Output}
\SetKwComment{Comment}{/*}{ */}
\let\oldnl\nl% Store \nl in \oldnl
\newcommand{\nonl}{\renewcommand{\nl}{\let\nl\oldnl}}

\theoremstyle{plain}

\theoremstyle{definition}

\theoremstyle{remark}

\usepackage[a4paper, top=1in, bottom=0.8in, left=0.7in, right=0.7in]{geometry}

\usepackage{graphicx}
\usepackage{epstopdf}% after graphicx
\usepackage{enumerate}
\usepackage{cite}
\usepackage{url} % not crucial - just used below for the URL
\usepackage{xcolor}
\usepackage{hyperref}
\usepackage{mathtools}
\usepackage{amsthm}
\usepackage{amsfonts}
\usepackage{amssymb}
\usepackage{amstext}
\usepackage{bbm}
\usepackage{enumitem}
\usepackage{hyperref}
\usepackage{float}
\usepackage{mathrsfs}

\usepackage{booktabs}
\usepackage{multirow}
\onecolumn

\begin{document}

\title{A Novel Methodology in Credit Spread Prediction Based on Ensemble Learning and Feature Selection}
\author{Yu~Shao,
%        ~\IEEEmembership{Member,~IEEE,}
        Jiawen~Bai,
%        ~\IEEEmembership{Fellow,~OSA,}
        Yingze~Hou,
        Xiaan~Zhou,
        and~Qihao~Pan
%        ~\IEEEmembership{Life~Fellow,~IEEE}% <-this % stops a space

\thanks{Y. Shao is with the Department of Mathematics and Statistics, Boston University, Boston, MA, 02215 USA e-mail: yshao19@bu.edu.}
\thanks{Y. Hou is with the Department of Industrial Engineering, University of Pittsburgh, Pittsburgh, PA, 15260 USA e-mail: yih81@pitt.edu.}
\thanks{X. Zhou is with Mays Business School, Texas A\&M University, College Station, TX, 77843 USA e-mail: xzhou520@tamu.edu.}% <-this % stops a space
% \thanks{Manuscript received April 19, 2005; revised August 26, 2015.}
}

% make the title area
\maketitle

\begin{abstract}
% The text of your abstract.  100 or fewer words.
The credit spread is a key indicator in bond investments, offering valuable insights for fixed-income investors to devise effective trading strategies. This study proposes a novel credit spread forecasting model leveraging ensemble learning techniques. To enhance predictive accuracy, a feature selection method based on mutual information is incorporated. Empirical results demonstrate that the proposed methodology delivers superior accuracy in credit spread predictions. Additionally, we present a forecast of future credit spread trends using current data, providing actionable insights for investment decision-making.
\end{abstract}


\IEEEpeerreviewmaketitle

\section{Introduction}

Credit spread has long been a critical focus for investors, particularly in the context of investment-grade corporate bonds, which have garnered even greater attention. Predicting the direction and magnitude of credit spread changes has been an enduring research topic.
The credit spread represents the yield difference between a risk-free security (typically a U.S. Treasury bond) and a risky security. A wider credit spread generally reflects a lower-quality bond, signaling a higher likelihood of issuer default. Consequently, credit spreads vary across securities depending on their credit ratings.
Beyond credit ratings, numerous other factors influence changes in credit spreads. Developing a practical and robust model to predict these changes, specifically for investment-grade corporate bonds, is therefore of significant value to investors.

Numerous studies have explored macroeconomic factors as determinants of credit spread changes, with research generally focusing on either positive or negative relationships.
Krueger and Kenneth (2003) employed a rational Bayesian model and regression analysis to demonstrate a positive relationship between unexpected increases in employment and benchmark Treasury rates \cite{krueger2003markets}. Wu and Zhang (2008) proposed an internally consistent approach to quantify the linkages between market prices of systematic macroeconomic risks and the term structure of credit spreads \cite{wu2008no}. Similarly, Davies (2008) provided evidence through hypothesis testing on the adverse effects of high inflation on the performance of corporate bonds \cite{davies2008postwar}.
On the other hand, studies on negative relationships have yielded equally significant findings. Gertler (1991) established a connection between credit spreads and GNP growth using a simple reduced-form test \cite{gertler1990interest}. Tang and Yan (2010) found that credit spreads tend to narrow with increases in GDP growth and the Consumer Confidence Index (CCI) based on their empirical analysis \cite{tang2010market}. Additionally, Collin (2001) conducted a comprehensive study incorporating a range of macro-level factors in a regression model to analyze both positive and negative influences on credit spreads \cite{collin2001determinants}.

Despite extensive research, significant uncertainties remain regarding the determinants of credit spread changes. One notable limitation is that existing factors and regression models explain only about a quarter of the observed variations \cite{collin2001determinants}. Furthermore, much of the analysis regarding these determinants is confined to qualitative insights, leaving a gap in quantitative prediction capabilities.

The primary contributions of this research are as follows:
\begin{enumerate}
    \item We introduced ensemble learning methods to address the limitations of existing credit spread forecasting approaches. Ensemble learning, recognized for its superior performance in financial forecasting compared to traditional machine learning algorithms \cite{ballings2015evaluating}, is applied to predict credit spread changes, an area that has received little attention thus far.

    \item Our methodology integrates ensemble learning with feature selection, providing valuable insights into the factors driving credit spread changes and offering actionable guidance for transaction decisions.

    \item Empirical results validate the effectiveness of our approach, demonstrating high accuracy in forecasting credit spread changes and establishing it as a promising solution to this persistent challenge.
    
\end{enumerate}

The rest of the paper is organized as follows. In Section \ref{sec: Methodology}, we present the methodology, covering feature description (which outlines the potential factors influencing credit spread changes), feature selection (which involves removing insignificant features), and the forecast model (which explains how the machine learning algorithm is applied). In Section \ref{sec: Empirical Analysis}, we conduct an empirical study and analyze the results to demonstrate the robustness of our model. Finally, in Section \ref{sec: Conclusion}, we provide our conclusion.

\section{Methodology}\label{sec: Methodology}

We propose a credit spread forecast model based on the ensemble learning method. First, we collect features from six key aspects closely related to credit spread changes, which form our raw feature set. Since the raw feature set may contain irrelevant information that introduces unnecessary noise into our prediction, we introduce a feature selection filter based on mutual information. This process not only helps in reducing noise but also highlights the contribution of each feature to credit spread prediction. The resulting filtered feature set, combined with the historical behavior of the credit spread, serves as the input for our prediction model.

\begin{figure}[h!]
\centering
\includegraphics[scale=0.36]{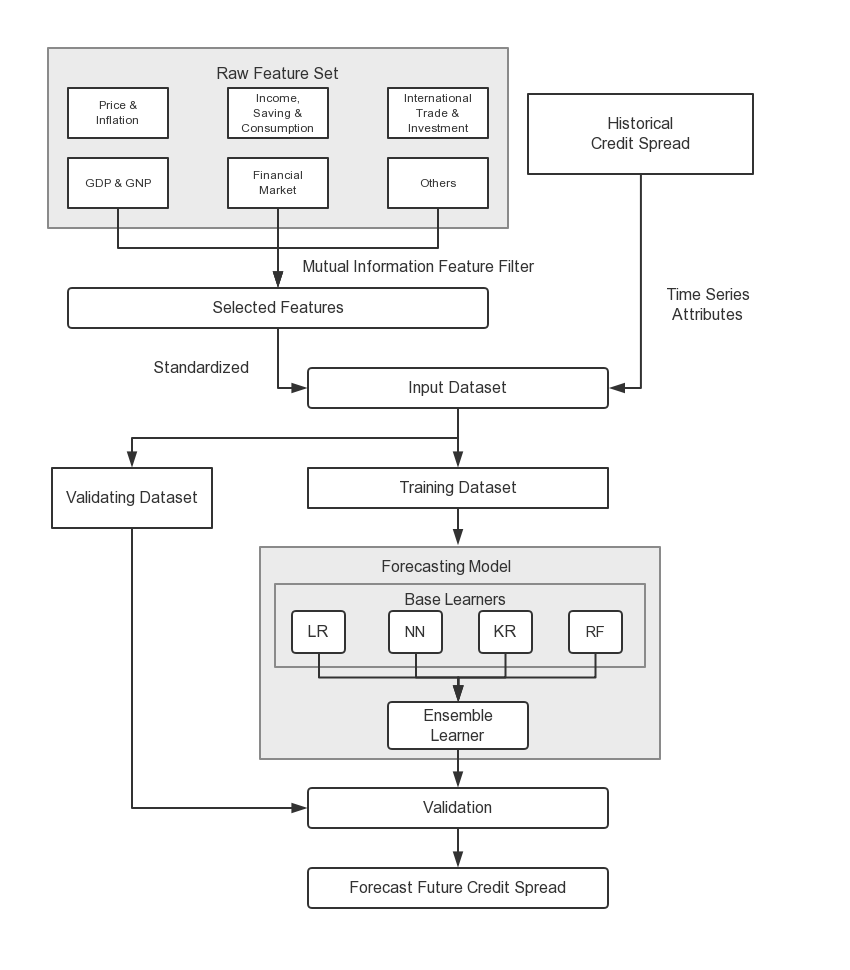}
\caption{The Framework of Credit Spread Forecasting Model}
\label{fig:methodology}
\end{figure}

Next, we apply several machine learning techniques and ensemble learning methods to forecast future credit spread changes. An overview of our methodology framework is presented in Figure \ref{fig:methodology}.
\vspace{-15pt}
\vspace{-10pt}

\subsection{Feature Description}

We identified 34 features that influence the credit spread of corporate bonds. For the sake of analysis, we categorize these features into six groups, each with an explanation.

To begin, let's review the features used to train our model. Gertler (1991) identified GNP as a factor influencing credit spread changes \cite{gertler1990interest}. A study by Collin (2001) summarized several financial market indicators that impact credit spreads, including the change in yield on 10-year Treasury bonds, the change in the 10-year minus 2-year Treasury yield spread, the change in implied volatility of the S\&P 500, and the return on the S\&P 500 \cite{collin2001determinants}. Christiansen (2002) analyzed the relationship between macroeconomic indices, such as the Producer Price Index (PPI), and changes in credit spreads \cite{christiansen2002credit}.
Krueger and Kenneth (2003) explored the relationship between the benchmark Treasury interest rate, a key component of credit spreads, and the unemployment rate \cite{krueger2003markets}. Wu and Zhang (2008) demonstrated the link between GDP growth rate and changes in credit spreads \cite{wu2008no}. Davies (2008) examined the impact of inflation risk, as reflected by the Consumer Price Index (CPI), on credit spreads \cite{davies2008postwar}.
Cúrdia (2010) studied the relationship between credit spreads and monetary policy, incorporating indices such as M1 and M2 money supply, government purchases, and government revenue \cite{curdia2010credit}. The work of Gilchrist (2010) confirmed the impact of significant macroeconomic fluctuations, such as industrial production, personal disposable income, and personal income, on credit spreads \cite{gilchrist2012credit}. Tang and Yan (2010) analyzed the effects of GDP growth rate and the Consumer Confidence Index (CCI) on credit spreads \cite{tang2010market}. Tsai (2010) used international trade indices to analyze the trends in credit spreads \cite{tsai2010combining}.
Table \ref{FeaturesTable1} lists the financial and economic indices referenced in some of these studies.

\begin{table}[h!]
\centering
\begin{tabular}{lp{11cm}}
\toprule
Author (Year)& Features\\
\midrule
Gertler (1991) & GNP\\
Collin-Dufresne et al. (2001) & Change in yield on 10-year Treasury; Change in 10-year minus 2-year Treasury yields; Change in implied volatility of S\&P 500; Return on S\&P 500     \\
Christiansen (2002) & Produce Price Index (PPI) \\
Krueger \& Kenneth (2003) & Unemployment rate \\
Davies (2008) & Inflation risk (CPI) \\
Wu \& Zhang (2008) & GDP; Real GDP \\
Cúrdia (2009) & M1 and M2 money stock; Government purchase;\\ & Government revenue \\
Gilchrist (2010) & Industrial Production; Personal disposable expenditure;\\ & Personal income \\
Tang \& Yan (2010) & Consumer Confidence Index (CCI) \\
Tsai (2010) & Net export volume; Export price index; Import price index; Total trading volume \\
\bottomrule
\end{tabular}
\caption{Features Associated with Credit Spread}
\label{FeaturesTable1}
\end{table}

We have gathered data for 34 features, primarily sourced from the St. Louis Fed data repository, with additional data obtained from the IMF Data and other databases. Based on the classification provided by the U.S. Bureau of Economic Analysis, we group the 34 features into six categories: price and inflation, income, saving and consumption, international trade and investment, GDP and GNP, financial market, and other factors. These categories are detailed in Table \ref{FeaturesTable2}. In addition, we also include the differences of each feature as separate variables during the feature selection process. Below is a brief overview of the six feature categories:

\begin{enumerate} 
    \item The "Price and Inflation" category quantifies the general price level and cost of living in the economy and measures the inflation rate. 
    \item The "Income, Saving, and Consumption" category focuses on personal income and consumption, providing insights into the economic situation in the U.S. 
    \item The "International Trade and Investment" category tracks changes in prices that firms and countries receive for products, reflecting global economic conditions. 
    \item The "GDP and GNP" category is self-explanatory, representing the nation's total economic output. 
    \item The "Financial Market" category includes indices commonly used as indicators of financial market performance and activity. 
    \item The "Other Factors" category comprises additional features that we and some prior studies believe have practical significance and the potential to represent economic trends in various aspects, though they are not substantial enough to warrant their own category. 
\end{enumerate}

\begin{table}[h!]
\centering
\begin{tabular}{lp{13cm}}
\toprule
Category & Features Name\\
\midrule
Price and Inflation 
    & CPI, Producer Price Index (PPI), Total Wholesale Trade Industries (WPI)\\
    & Consumer Confidence Index (CCI), Producer Price Index (PPI)\\
    & GDP Deflator ($GDP_D$)\\
\midrule
Income, Saving & Personal Consumption Expenditures (PCE), Personal Saving Rate (PSR)\\
and Consumption	 
    & Disposable Personal Income (DPI)\\
\midrule
International Trade
& Export Price Index (EPI), Export volume index (EVI), Export rate (ER)\\
and Investment 
	& Import Price Index (IPI), Import volume index (IVI), Import rate (IR)\\
	& Price Index of Machinery product export\\
	& Quantity Index of Machinery product export\\
	& Electric product export order, Electric machinery product export order\\
	& Information and communication product export order\\
	& Net Exports of Goods and Services (NEGS), Total trading volume Index (TTVIND)\\
\midrule
GDP and GNP 
    & GDP, Real GDP (R\_GDP), GNP, Percentage change in GDP ($\Delta$GDP)\\
    & Percentage change in Real GDP ($\Delta$RGDP), Percentage change in GNP ($\Delta$GNP), Industrial production Index (INDPRO)\\
\midrule
Financial Market & 10-Year Treasury Constant Maturity Minus 2-Year Treasury Constant Maturity (T10Y2YM), 10-Year Treasury Constant Maturity Rate (GS10)\\
	& VIX index (VIX), S\&P 500 return (S\&P500), Discount Rate (DR)\\
	& M1, M2, USD Index (USDIND)\\
\midrule
Other Factors 
    & Leading Index (LIND), Lagging Index (LIND2), Unemployment rate (UR)\\
\bottomrule
\end{tabular}
\caption{Category of Features}
\label{FeaturesTable2}
\end{table}

\subsection{Feature Selection}
In our feature selection method, the criterion for selecting features is the amount of information a feature can contribute to the forecast system. The more information a feature provides, the more crucial it is for the prediction process. For each feature, the amount of information changes depending on whether or not it is included in the forecast system. The difference in the information amount between the system with and without the feature represents the information gain brought by the feature. This information gain guides us in making better forecasting decisions. To measure the amount of information in the context of forecasting, we use differential entropy.

Given a random variable $X$ with probability density function $f(x)$, the differential entropy $h(X)$ is defined as
\begin{equation}
    h(X) = - \int f(x) \log f(x) dx. \nonumber
\end{equation}
For two random variables $X$ and $Y$, suppose they have a joint pdf $f(x,y)$, their conditional differential entropy $h(X|Y)$ is defined as
\begin{equation}
    h(X|Y) = - \int f(x,y) \log f(x|y) dx dy \nonumber
\end{equation}
and the mutual information between $X$ and $Y$ is given as
\begin{equation}
    I(X ; Y) = h(X) - h(X|Y). \nonumber
\end{equation}

Notice that $I(X;Y) \geq 0$ and $h(X|Y) \leq h(X)$ where the equalities hold if and only if $X$ and $Y$ are independent. If we regard the future credit spread as the random variable $X$, then the entropy $h(X)$ measures the randomness of $X$. The inclusion of any feature $Y$ will reduce the randomness of $X$, thereby improving our ability to predict the credit spread. The higher the mutual information between $X$ and $Y$, the more information $Y$ provides to the prediction. Therefore, we aim to calculate the mutual information between each feature $Y$ and the credit spread $X$ using historical data, and select the features with the highest mutual information values as the inputs for our prediction model.

\begin{figure}[h!]
\centering
\includegraphics[scale=0.75]{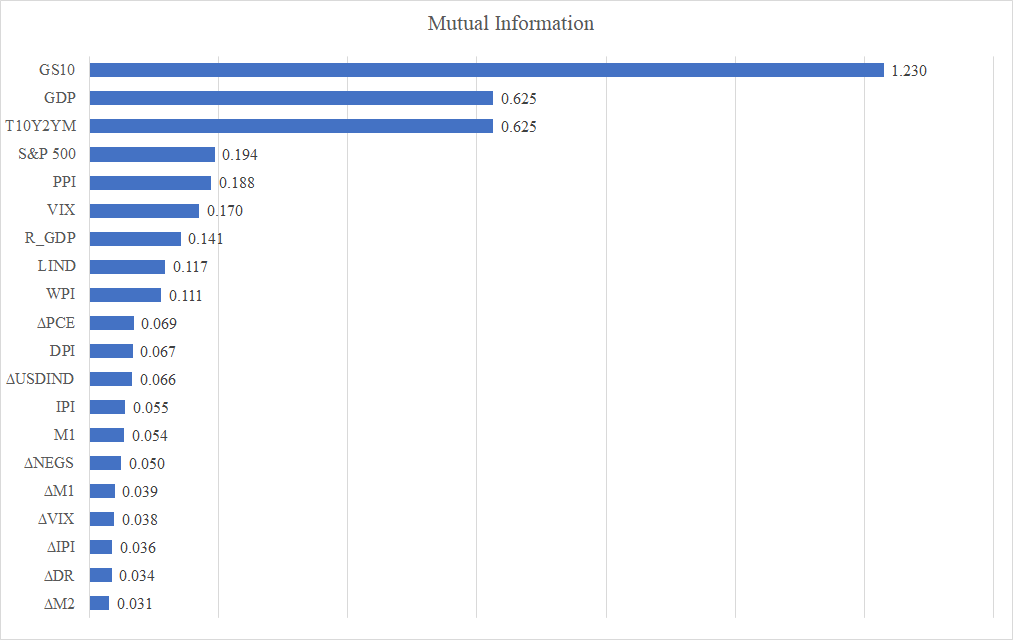}
\caption{Mutual Information of Selected Features}
\label{fig:selectedfeature}
\end{figure}

Figure \ref{fig:selectedfeature} shows the features with the highest mutual information. We select 20 features from the raw feature set to construct our prediction model in order to avoid disturbance from insignificant features. Since mutual information reflects the significance of features to the credit spread, features with higher mutual information are more helpful for predicting the future spread. Therefore, from Figure \ref{fig:selectedfeature}, the features with the highest influence on the credit spread are the 10-Year Treasury Constant Maturity Rate, GDP, 10-Year Treasury Constant Maturity Minus 2-Year Treasury Constant Maturity, S\&P 500 index, PPI, and VIX index, which correspond to our analysis in the previous section.

\subsection{Prediction Model}
In this section, we will set up a two-layer prediction model. The forecasting model consists of several base learners including Multi-layer Perceptron regressor (MLP), random forest regressor, k-nearest neighbors regressor (K-NN), and an ensemble learner as the second layer.

In supervised learning algorithms of machine learning, our goal is to develop a stable model that performs well in all aspects. However, the actual situation is often not so ideal; sometimes we can only obtain multiple weak supervised models that perform better in certain aspects. Ensemble learning combines these weak models to create a better and more comprehensive model. When spread price changes dramatically, ensemble learning prevents the original model from being affected by outlier values. There are several ensemble learning techniques available, and stacking will be used in our work. The algorithm in Table \ref{AlgorithmTable1} summarizes stacking.

\begin{table}[h!]
\centering
\begin{tabular}{lp{11cm}}

\toprule
\textbf{Algorithm:} Stacking \\
\midrule
1: \textbf{Input:} training data $D=\{x_i, y_i\}_{i=1}^m $ \\
2: \textbf{Ouput:} ensemble classifier $H$ \\
3: \textit{Step 1: learn base-level regressor} \\
4: \textbf{for} $t = 1$ to $T$ \textbf{do} \\ 
5:\ \ \ learn $h_t$ based on $D$ \\
6: \textbf{end for} \\
7: \textit{Step 2: construct new data set of predictions} \\
8: \textbf{for} $t = 1$ to $T$ \textbf{do} \\ 
9:\ \ \ $D_h = \{x_i^{'}, y_i\}$ where $x_i^{'} = \{h_1(x_i), ..., h_T(x_i)\}$ \\
10: \textbf{end for} \\
11: \textit{Step 3: learn a meta-regressor} \\
12: learn $H$ based on $D_h$ \\
13: return $H$ \\
\bottomrule
\end{tabular}
\caption{Stacking Algorithm}
\label{AlgorithmTable1}
\end{table}

In the first layer of our model, we should get predictions from three base learners. These predictions will be passed as features to the ensemble learner and will be trained in the second layer. Specifically, our predicting process is as follows:

\begin{enumerate}
    \item The spread price itself reveals much information about the future movement of the spread price. In addition to those features we have mentioned in Section II, we include the average spread over the past few months. While this limits the length of time we can predict, it will correct the future predicted value by updating the spread price. The optimal duration for the average spread price is calculated by minimizing the mean square error. Next, we will conduct PCA whitening with the recent average spread price and the other financial indicators. The goal here is to reduce the correlation between all features, as financial data are usually highly correlated. Then, the whitened data will be passed to Multi-layer Perceptron regressor (MLP), Random Forest regressor, and K-NN regression respectively. The predictions from these regressors will be used to train the Kernel Ridge regressor, which is the second layer of our model.

    \item These three base learners are chosen because of their excellent predictive effect on the historical spread price. MLP is a class of feedforward artificial neural network, and it allows us to solve problems stochastically, which makes it a good regression method for the spread price. Random Forest regressor is chosen for its popularity among research scientists and its high accuracy. K-NN works perfectly around local values. The limitations of each method drive us to set up a second layer (Kernel Ridge) to balance the predictions as well as reduce the noise.
\end{enumerate}

\section{Empirical Analysis} \label{sec: Empirical Analysis}

The empirical data in our work comes from websites of the Bureau of Economic Analysis and the Federal Reserve. We collect 120 pieces of monthly historical data lasting from Jan 2008 to Dec 2017 to construct our experimental data set. Each piece of data contains the features mentioned in Table \ref{FeaturesTable2} as the input parameters and the monthly credit spread as the target of the learning model. We divide our data set into two parts, where the first 70\% of the data construct the training set for the modeling process, and the latter 30\% construct the testing set for validation. For each prediction model, three indicators, the Mean Absolute Error (MAE), Mean Squared Error (MSE), and $R^2$ score, are employed as the criteria for the performance evaluation of models. Given $n$ pairs of actual values $y_i$ and the predicted values $\hat{y}_i$, the indicators are described as follows:
\begin{align*}
     MAE = \frac{1}{n} \sum_{i=1}^n \left| y_i - \hat{y}_i \right|, \quad\quad
    MSE = \frac{1}{n} \sum_{i=1}^n \left( y_i - \hat{y}_i \right)^2, \quad\quad
    R^2 = 1- \frac{\sum_{i=1}^n (y_i - \hat{y}_i)^2}{\sum_{i=1}^n (y_i-\bar{y})^2}.
\end{align*}

\subsection{Experimental Result}

We compare the prediction results of four base learners and the ensemble learner. For each learner, we consider cases of using the raw features and using the selected feature set. The results of the training and testing performances are displayed in Table \ref{Performance}.

As we can see from Table \ref{Performance}, the best prediction on both the training set and testing set appears in the stacking model. Notice that in most algorithms, although feature selection doesn't improve the results on the training set, the results on the testing set are better. This is because the overfitting problem on the original data set is avoided by feature selection, and the anti-noise ability of the model has also been improved. Moreover, besides the Stacking algorithm, all methods have their own limitations. For linear regression, although it fits the training set well, the predictions are even worse than using the mean prices on the testing data set. For K-NN regression, the drawback lies in its sensitivity to the local structure of the dataset. For kernel ridge, the results are not bad, but there is still room for improvement. The random forest also has decent predictions, but it is subject to randomness and is not robust on the data set. 

\begin{table}[h!]
\centering
\begin{tabular}{|c|c|c|c|c|c|c|c|}
\hline
\multirow{2}{*}{Learning Method} & \multirow{2}{*}{Feature Selection} & \multicolumn{3}{c|}{Training Set Result} & \multicolumn{3}{c|}{Testing Set Result} \\ \cline{3-8} 
 &  & MAE & MSE & $R^2$ Score & MAE & MSE & $R^2$ Score \\ \hline
\multirow{2}{*}{Linear Regression} & No & 0.082 & 0.010 & 0.981 & 0.528 & 1.378 & -0.787 \\ \cline{2-8} 
 & Yes & 0.086 & 0.011 & 0.978 & 0.463 & 1.163 & -0.508 \\ \hline
\multirow{2}{*}{K-NN regression} & No & 0.119 & 0.127 & 0.752 & 0.180 & 0.120 & 0.844 \\ \cline{2-8} 
 & Yes & 0.119 & 0.127 & 0.752 & 0.180 & 0.120 & 0.844 \\ \hline
\multirow{2}{*}{Kernel Ridge} & No & 0.119 & 0.027 & 0.948 & 0.235 & 0.111 & 0.856 \\ \cline{2-8} 
 & Yes & 0.121 & 0.028 & 0.945 & 0.231 & 0.107 & 0.861 \\ \hline
\multirow{2}{*}{Random Forest} & No & 0.104 & 0.064 & 0.874 & 0.238 & 0.165 & 0.786 \\ \cline{2-8} 
 & Yes & 0.098 & 0.064 & 0.875 & 0.204 & 0.106 & 0.862 \\ \hline
\multirow{2}{*}{Stacking} & No & 0.003 & 1e-4 & 0.999 & 0.168 & 0.071 & 0.908 \\ \cline{2-8} 
 & Yes & 0.004 & 1e-4 & 0.999 & 0.155 & 0.062 & 0.920 \\ \hline
\end{tabular}
\caption{The Performance of Credit Spread Prediction Models}
\label{Performance}
\end{table}

\begin{figure}[h!]
\centering
\begin{minipage}{0.4\textwidth}
    \centering
    \includegraphics[width=\textwidth]{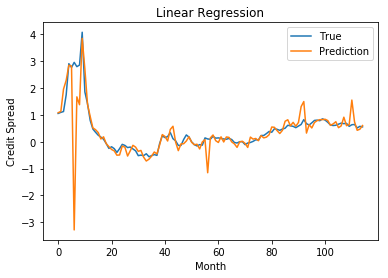}
    \caption*{(a) Linear Regression}
\end{minipage}
% \hfill
\begin{minipage}{0.4\textwidth}
    \centering
    \includegraphics[width=\textwidth]{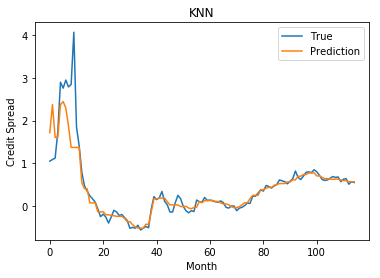}
    \caption*{(b) K-NN Regression}
\end{minipage}

\vspace{1em}

\begin{minipage}{0.4\textwidth}
    \centering
    \includegraphics[width=\textwidth]{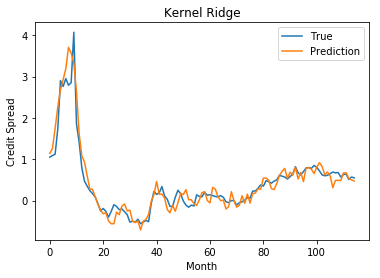}
    \caption*{(c) Kernel Ridge Regression}
\end{minipage}
% \hfill
\begin{minipage}{0.4\textwidth}
    \centering
    \includegraphics[width=\textwidth]{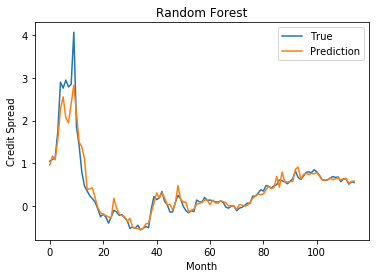}
    \caption*{(d) Random Forest Regression}
\end{minipage}

\caption{Prediction of Base Learners}
\label{fig:base-learners}
\end{figure}

Also, we plot a simulated price prediction and provide a short analysis of each base learner in Figure \ref{fig:base-learners}.
In Figure \ref{fig:base-learners}(a), the linear regression could be seen as a baseline in this work, but it does not produce a convincing result from the plot; the predictions fluctuate around the true prices. In Figure \ref{fig:base-learners}(b), K-NN regression smoothes the historical spread and provides a better estimate in the period of relatively stable prices. However, for the estimation of the spread price with high volatility during the financial crisis period, the estimation is not satisfying due to the sensitivity to the local structure. In Figure \ref{fig:base-learners}(c), kernel ridge is an improved version of the linear regression method here, where we add a penalty term to the features to avoid the overfitting problem of regression. Therefore, we have a better output than the linear regression. In Figure \ref{fig:base-learners}(d), the random forest is actually an ensemble learning method, and it indeed provides a decent result, which gives us confidence to believe that the stacking algorithm should work in this context.

\begin{figure}[h!]
\centering
\begin{minipage}[b]{0.4\textwidth}
    \centering
    \includegraphics[width=\textwidth]{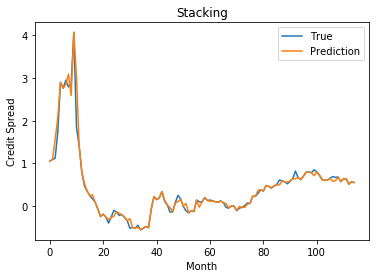}
    \caption{Stacking}
    \label{fig:Stacking}
\end{minipage}
\hspace{0.02\textwidth} % Adjust the space between the two figures
\begin{minipage}[b]{0.4\textwidth}
    \centering
    \includegraphics[width=\textwidth]{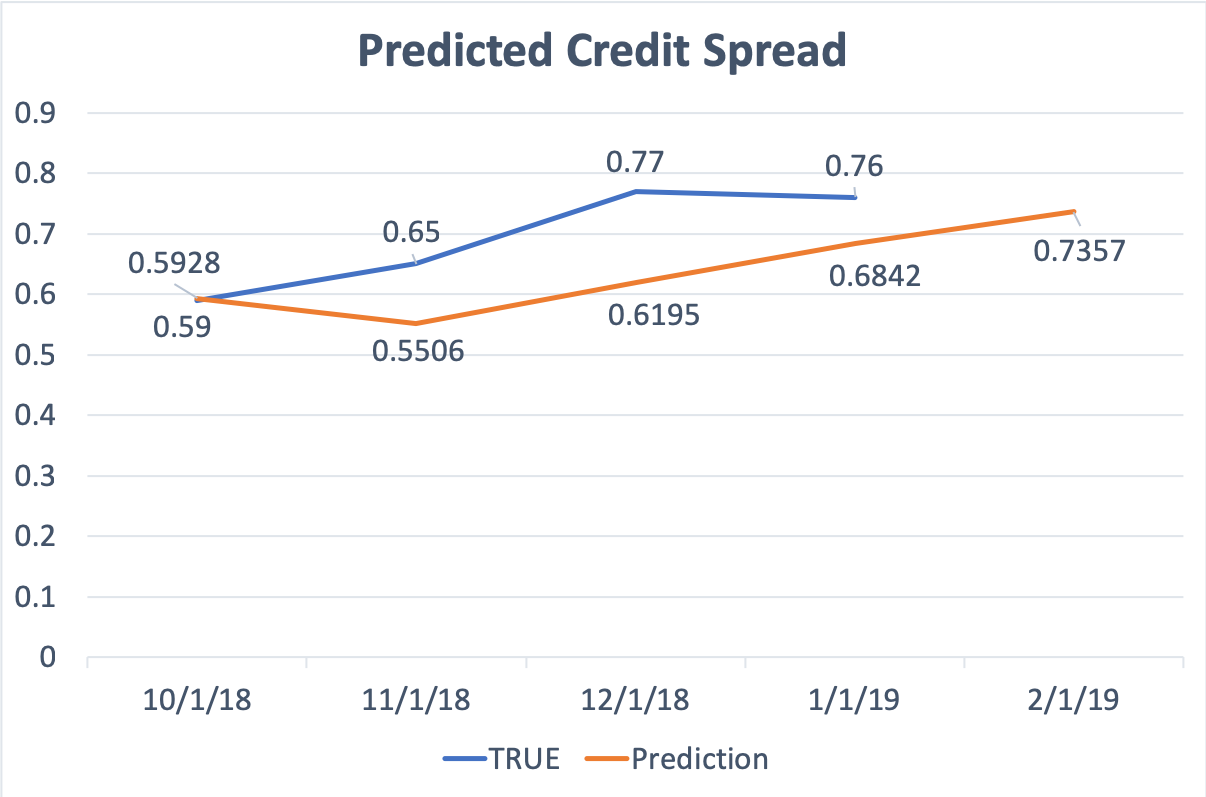}
    \caption{Prediction of Credit Spread}
    \label{fig:Prediction}
\end{minipage}
\end{figure}

In Figure \ref{fig:Stacking}, the stacking method combines the output from other predictors, providing us with an algorithm robust enough to handle most situations, and it also yields the best outcome. As we can see from the plot, not only does it perform well on the training set, but it also works very well on the testing set.

\subsection{Prediction}

From the previous discussion, we know that the ensemble learning model performs the best when feature selection techniques are included. Therefore, given the current data, we use the stacking method with the selected feature set to predict the future behavior of the credit spread. With the data from 2018, we predict the credit spread for the first quarter of 2019, as shown in Figure \ref{fig:Prediction}. The actual credit spreads for the previous four months validate our prediction, with errors less than $15$ bps. We also provide a forecast for the credit spread of February 2019 to be $73$ bps. In general, the predicted value is less than the true value, which is due to some null values in some features affecting the accuracy. However, the predictions regarding the direction and magnitude of the credit spread movement are quite accurate.

\section{Conclusion} \label{sec: Conclusion}
This paper illustrates a novel prediction method for the monthly credit spread based on the ensemble learning method. 
We innovatively include the feature selection method using mutual information in forecasting the credit spread. 
The empirical results show that the ensemble learning method performs better than the traditional machine learning methods. 
We predict the credit spread in February 2019 to be $73$ bps. 
Moreover, feature selection before prediction not only explains the rationality of the features but also improves the accuracy and robustness of the prediction results. 

However, our work has some limitations. First, although mutual information shows the significance of the selected features, it is difficult to identify the cross relationships within different features, which requires further feature engineering methods. Second, we need cross-sectional data for the prediction, meaning we must have access to the most recent monthly data to forecast the credit spread for the next month. This limitation of our prediction model restricts the time frame for future predictions. Future work could incorporate more efficient time series analysis techniques into our model to improve performance and extend the forecasting period.

\bibliographystyle{plain}
\bibliography{references}
\end{document}